\documentclass[a4paper, 12 pt]{article}

\usepackage[T2A]{fontenc}
\usepackage[cp1251]{inputenc}
\usepackage[english]{babel}
\usepackage[tbtags]{amsmath}
\usepackage{amsfonts,amssymb}

\usepackage{mathrsfs}

\usepackage{geometry} 
\begin{document}
\title{Riordan arrays, Chebyshev polynomials, Fibonacci bases}
 \author{E. Burlachenko}
 \date{}

 \maketitle
\begin{abstract}
Chebyshev polynomials and their modifications are attributes of various fields of mathematics. In particular, they are generating functions of the rows elements of certain Riordan matrices. In paper, we give a selection of some characteristic situations in which such matrices are involved. Using the columns and rows of these matrices, we will build  the bases of the space of formal power series and the space of polynomials, the properties of which allow us to call them "Fibonacci bases".
\end{abstract}
\section{Introduction}
Subject of our study are the transformations in space of formal power series and the corresponding matrices. We will associate rows and columns of matrices with the generating functions of their elements, i.e. with the formal power series. $n$th coefficient of the series $a\left( x \right)$, $n$th row and $n$th column of the matrix $A$ will be denoted  respectively by
$$\left[ {{x}^{n}} \right]a\left( x \right),\qquad    \left[ n,\to  \right]A,\qquad   A{{x}^{n}}.$$
Matrix $A$, transformation corresponding to the matrix $A$ and the basis of vector space formed by the sequence of columns of the matrix $A$ will be denoted by the same symbol. I.e. we will say: matrix $A$, transformation $A$, basis $A$.

Infinite lower triangular matrix $\left( f\left( x \right),g\left( x \right) \right)$ $n$th column of which has the generating function $f\left( x \right){{g}^{n}}\left( x \right)$, ${{g}_{0}}=0$, is called Riordan array, or Riordan matrix [1] – [3]. It is the product of two matrices that correspond to the operators of multiplication and composition of series:
$$\left( f\left( x \right),g\left( x \right) \right)=\left( f\left( x \right),x \right)\left( 1,g\left( x \right) \right),$$
$$\left( f\left( x \right),x \right)a\left( x \right)=f\left( x \right)a\left( x \right), \qquad\left( 1,g\left( x \right) \right)a\left( x \right)=a\left( g\left( x \right) \right),$$
$$\left( f\left( x \right),g\left( x \right) \right)\left( b\left( x \right),a\left( x \right) \right)=\left( f\left( x \right)b\left( g\left( x \right) \right),a\left( g\left( x \right) \right) \right).$$
Transposed Riordan matrices can be considered as matrices of the operators acting in the space of polynomials. Similar operators corresponding to the transposed exponential Riordan matrices (i.e. matrices $\left| {{e}^{x}} \right|{{\left( f\left( x \right),g\left( x \right) \right)}^{T}}{{\left| {{e}^{x}} \right|}^{-1}}$, where $\left| {{e}^{x}} \right|$ is the diagonal matrix: $\left| {{e}^{x}} \right|{{x}^{n}}={{{x}^{n}}}/{n!}\;$) are considered in the umbral calculus [4].

Riordan matrix inverse to itself is called the Riordan involution [5], [6]. It can be represented in the form $RM$, where $R$ is the some Riordan matrix,
$$M=\left( 1,-x \right)=\left( \begin{matrix}
   1 & 0 & 0 & 0 & \cdots   \\
   0 & -1 & 0 & 0 & \cdots   \\
   0 & 0 & 1 & 0 & \cdots   \\
   0 & 0 & 0 & -1 & \cdots   \\
   \vdots  & \vdots  & \vdots  & \vdots  & \ddots   \\
\end{matrix} \right).$$
Matrix $R$ is called the Riordan pseudo-involution. We will use this terminology also for transposed Riordan  matrices. Matrix
$$P=\left( \frac{1}{1-x},\frac{x}{1-x} \right)=\left( \begin{matrix}
   1 & 0 & 0 & 0 & \cdots   \\
   1 & 1 & 0 & 0 & \cdots   \\
   1 & 2 & 1 & 0 & \cdots   \\
   1 & 3 & 3 & 1 & \cdots   \\
   \vdots  & \vdots  & \vdots  & \vdots  & \ddots   \\
\end{matrix} \right)$$
is called Pascal matrix. Power of the Pascal matrix defined by ${{P}^{\varphi }}=\left( \frac{1}{1-\varphi x},\frac{x}{1-\varphi x} \right)$. Transformation ${{P}^{\varphi }}$ acts in the space of formal power series and is called the generalized Euler transformation:
$${{P}^{\varphi }}a\left( x \right)=\frac{1}{1-\varphi x}a\left( \frac{x}{1-\varphi x} \right);$$
transformation${{\left( {{P}^{\varphi }} \right)}^{T}}$ acts in the space of polynomials and is called the shift operator:
 $${{\left( {{P}^{\varphi }} \right)}^{T}}c\left( x \right)=c\left( x+\varphi  \right).$$ 
Transformation ${{\left( {{P}^{\varphi }}M \right)}^{T}}=M{{\left( {{P}^{\varphi }} \right)}^{T}}$ is an involution. Each polynomial $c\left( x \right)$ is the sum of polynomials ${{c}_{1}}\left( x \right)$, ${{c}_{2}}\left( x \right)$ such that $M{{\left( {{P}^{\varphi }} \right)}^{T}}{{c}_{1}}\left( x \right)={{c}_{1}}\left( x \right)$, $M{{\left( {{P}^{\varphi }} \right)}^{T}}{{c}_{2}}\left( x \right)=-{{c}_{2}}\left( x \right)$. Hence, ${{\left( {{P}^{\varphi }} \right)}^{T}}{{c}_{1}}\left( x \right)={{c}_{1}}\left( -x \right)$, ${{\left( {{P}^{\varphi }} \right)}^{T}}{{c}_{2}}\left( x \right)=-{{c}_{2}}\left( -x \right)$. Respectively, since transformation $M{{P}^{\varphi }}={{P}^{-\varphi }}M$ is an involution, each series $a\left( x \right)$ is the sum of series ${{a}_{1}}\left( x \right)$, ${{a}_{2}}\left( x \right)$ such that ${{P}^{\varphi }}{{a}_{1}}\left( x \right)={{a}_{1}}\left( -x \right)$, ${{P}^{\varphi }}{{a}_{2}}\left( x \right)=-{{a}_{2}}\left( -x \right)$. Sequence of columns of the invertible infinite upper triangular matrix $B$ such that ${{\left( {{P}^{\varphi }} \right)}^{T}}B=MBM$ will be called pseudo-eigenbasis of the transformation ${{\left( {{P}^{\varphi }} \right)}^{T}}$; sequence of columns of the invertible infinite lower triangular matrix $A$ such that ${{P}^{\varphi }}A=MAM$ will be called pseudo-eigenbasis of the transformation ${{P}^{\varphi }}$ . Polynomials ${{c}_{1}}\left( x \right)$, ${{c}_{2}}\left( x \right)$ are the combinations respectively of even and odd columns of the matrix $B$; series ${{a}_{1}}\left( x \right)$, ${{a}_{2}}\left( x \right)$ are the combinations respectively of even and odd columns of the matrix $A$.

Eigenspaces of the transformation$PM$ are subject of papers [7] – [9]. In [10], [11] eigenspaces of the transformations $PM$ and ${{P}^{T}}M$ are considered on general terms; this point of view intersects with our observations set out in Section 4.

Modified Chebyshev polynomials of the first and second kind ${{C}_{n}}\left( x \right)=2{{T}_{n}}\left( {x}/{2}\; \right)$, ${{S}_{n}}\left( x \right)={{U}_{n}}\left( {x}/{2}\; \right)$, or
$${{C}_{n}}\left( x \right)=\prod\limits_{m=1}^{n}{\left( x-2\cos \frac{2m-1}{2n}\pi  \right)},\qquad  {{S}_{n}}\left( x \right)=\prod\limits_{m=1}^{n}{\left( x-2\cos \frac{m}{n+1}\pi  \right)},$$
${{C}_{0}}\left( x \right)=2$, ${{S}_{0}}\left( x \right)=1$, are associated with the Riordan matrices as follows:
 $${{C}_{n}}\left( x \right)=\left[ n,\to  \right]\left( \frac{1-{{x}^{2}}}{1+{{x}^{2}}},\frac{x}{1+{{x}^{2}}} \right),\quad n>0;\qquad  {{S}_{n}}\left( x \right)=\left[ n,\to  \right]\left( \frac{1}{1+{{x}^{2}}},\frac{x}{1+{{x}^{2}}} \right),$$ 
$$\left( \frac{1-{{x}^{2}}}{1+{{x}^{2}}},\frac{x}{1+{{x}^{2}}} \right)=\left( \begin{matrix}
   1 & 0 & 0 & 0 & 0 & 0 & 0 & \cdots   \\
   0 & 1 & 0 & 0 & 0 & 0 & 0 & \cdots   \\
   -2 & 0 & 1 & 0 & 0 & 0 & 0 & \cdots   \\
   0 & -3 & 0 & 1 & 0 & 0 & 0 & \cdots   \\
   2 & 0 & -4 & 0 & 1 & 0 & 0 & \cdots   \\
   0 & 5 & 0 & -5 & 0 & 1 & 0 & \cdots   \\
   -2 & 0 & 9 & 0 & -6 & 0 & 1 & \cdots   \\
   \vdots  & \vdots  & \vdots  & \vdots  & \vdots  & \vdots  & \vdots  & \ddots   \\
\end{matrix} \right),$$
$$\left( \frac{1}{1+{{x}^{2}}},\frac{x}{1+{{x}^{2}}} \right)=\text{ }\left( \begin{matrix}
   1 & 0 & 0 & 0 & 0 & 0 & 0 & \cdots   \\
   0 & 1 & 0 & 0 & 0 & 0 & 0 & \cdots   \\
   -1 & 0 & 1 & 0 & 0 & 0 & 0 & \cdots   \\
   0 & -2 & 0 & 1 & 0 & 0 & 0 & \cdots   \\
   1 & 0 & -3 & 0 & 1 & 0 & 0 & \cdots   \\
   0 & 3 & 0 & -4 & 0 & 1 & 0 & \cdots   \\
   -1 & 0 & 6 & 0 & -5 & 0 & 1 & \cdots   \\
   \vdots  & \vdots  & \vdots  & \vdots  & \vdots  & \vdots  & \vdots  & \ddots   \\
\end{matrix} \right).$$
Generalization of these polynomials are the Dickson polynomials:
$${{D}_{n}}\left( x,\beta  \right)=\left[ n,\to  \right]\left( \frac{1-\beta {{x}^{2}}}{1+\beta {{x}^{2}}},\frac{x}{1+\beta {{x}^{2}}} \right)=\prod\limits_{m=1}^{n}{\left( x-2\sqrt{\beta }\cos \frac{2m-1}{2n}\pi  \right)},\quad n>0,$$
$${{E}_{n}}\left( x,\beta  \right)=\left[ n,\to  \right]\left( \frac{1}{1+\beta {{x}^{2}}},\frac{x}{1+\beta {{x}^{2}}} \right)=\prod\limits_{m=1}^{n}{\left( x-2\sqrt{\beta }\cos \frac{m}{n+1}\pi  \right)},$$
${{D}_{0}}\left( x,\beta  \right)=2$, ${{E}_{0}}\left( x \right)=1$. Hence, $n$tn row of the matrix 
$$\left( \frac{1-\beta {{x}^{2}}}{1-\varphi x+\beta {{x}^{2}}},\frac{x}{1-\varphi x+\beta {{x}^{2}}} \right)=\left( \frac{1-\beta {{x}^{2}}}{1+\beta {{x}^{2}}},\frac{x}{1+\beta {{x}^{2}}} \right)\left( \frac{1}{1-\varphi x},\frac{x}{1-\varphi x} \right),$$
$n>0$, is the polynomial ${{D}_{n}}\left( x+\varphi ,\beta  \right)$, $n$tn row of the matrix 
$$\left( \frac{1}{1-\varphi x+\beta {{x}^{2}}},\frac{x}{1-\varphi x+\beta {{x}^{2}}} \right)=\left( \frac{1}{1+\beta {{x}^{2}}},\frac{x}{1+\beta {{x}^{2}}} \right)\left( \frac{1}{1-\varphi x},\frac{x}{1-\varphi x} \right)$$
is the polynomial ${{E}_{n}}\left( x+\varphi ,\beta  \right)$.

Chebyshev polynomials and their modifications in connection with Riordan matrices are considered in [12] – [15]. 

Generalized Fibonacci and Lucas sequences are associated with the Dickson polynomials and the generalized Euler transformation as follows. Denote
$$F_{n}^{\left( \varphi ,\beta  \right)}=\varphi F_{n-1}^{\left( \varphi ,\beta  \right)}-\beta F_{n-2}^{\left( \varphi ,\beta  \right)},  \qquad F_{0}^{\left( \varphi ,\beta  \right)}=0,\qquad F_{1}^{\left( \varphi ,\beta  \right)}=1,$$
$$L_{n}^{\left( \varphi ,\beta  \right)}=\varphi L_{n-1}^{\left( \varphi ,\beta  \right)}-\beta L_{n-2}^{\left( \varphi ,\beta  \right)},  \qquad L_{0}^{\left( \varphi ,\beta  \right)}=2, \qquad L_{1}^{\left( \varphi ,\beta  \right)}=\varphi .$$
Then
$$\sum\limits_{n=0}^{\infty }{L_{n}^{\left( \varphi ,\beta  \right)}}{{x}^{n}}=\frac{2-\varphi x}{1-\varphi x+\beta {{x}^{2}}},   \qquad\sum\limits_{n=0}^{\infty }{F_{n}^{\left( \varphi ,\beta  \right)}{{x}^{n}}}=\frac{x}{1-\varphi x+\beta {{x}^{2}}}.$$
Since
$$\left( \frac{1-\beta {{x}^{2}}}{1+\beta {{x}^{2}}},\frac{x}{1+\beta {{x}^{2}}} \right)\frac{1}{1-\varphi x}=\frac{1-\beta {{x}^{2}}}{1-\varphi x+\beta {{x}^{2}}}=\frac{2-\varphi x}{1-\varphi x+\beta {{x}^{2}}}-1,$$
$$\left( \frac{1}{1+\beta {{x}^{2}}},\frac{x}{1+\beta {{x}^{2}}} \right)\frac{1}{1-\varphi x}=\frac{1}{1-\varphi x+\beta {{x}^{2}}},$$
then
$$\frac{2-\varphi x}{1-\varphi x+\beta {{x}^{2}}}=\sum\limits_{n=0}^{\infty }{{{D}_{n}}}\left( \varphi ,\beta  \right){{x}^{n}},  \qquad\frac{x}{1-\varphi x+\beta {{x}^{2}}}=\sum\limits_{n=1}^{\infty }{{{E}_{n-1}}}\left( \varphi ,\beta  \right){{x}^{n}}.$$
Since
$${{P}^{-\varphi }}\frac{2-\varphi x}{1-\varphi x+\beta {{x}^{2}}}=\frac{2+\varphi x}{1+\varphi x+\beta {{x}^{2}}}, \qquad{{P}^{-\varphi }}\frac{x}{1-\varphi x+\beta {{x}^{2}}}=\frac{x}{1+\varphi x+\beta {{x}^{2}}},$$
transformation ${{P}^{-\varphi }}$ translates the series (geometric progression) 
$$\frac{1}{2}\left( \frac{2-\varphi x}{1-\varphi x+\beta {{x}^{2}}}+\frac{x\sqrt{{{\varphi }^{2}}-4\beta }}{1-\varphi x+\beta {{x}^{2}}} \right)={{\left( 1-\left( \frac{\varphi +\sqrt{{{\varphi }^{2}}-4\beta }}{2} \right)x \right)}^{-1}}={{\left( 1-\lambda x \right)}^{-1}}$$
in the geometric progression
$$\frac{1}{2}\left( \frac{2+\varphi x}{1+\varphi x+\beta {{x}^{2}}}+\frac{x\sqrt{{{\varphi }^{2}}-4\beta }}{1+\varphi x+\beta {{x}^{2}}} \right)={{\left( 1+\left( \frac{\varphi -\sqrt{{{\varphi }^{2}}-4\beta }}{2} \right)x \right)}^{-1}}={{\left( 1+\frac{\beta }{\lambda }x \right)}^{-1}}.$$

In the following sections of  paper, we present some observations related to the stated subject matter. In Sections 1 and 2, we consider two types of the transformations and corresponding Riordan matrices. Rows and columns of these matrices, firstly, are associated with the certain modifications of the Chebyshev polynomials, and secondly, they form the pseudo-eigenbases of the transformations ${{P}^{\varphi }}$, ${{\left( {{P}^{\varphi }} \right)}^{T}}$. In Section 4, we will build the pseudo-eigenbases of the transformations $P$ and ${{\left( {{P}^{-1}} \right)}^{T}}$, using the rows and columns of the matrices of certain “transformations of the second type ”. Some features of these bases allow us to call them "Fibonacci bases". In Section 5, we give a generalization of these bases. In all the cases considered, we will use two identities for the rows of Riordan matrices which constitute the content of the following theorem.\\
{\bfseries Theorem 1.} \emph{ If the series $a\left( x \right)$, ${{a}_{0}}=1$, $b\left( x \right)$, ${{b}_{0}}=1$, are bound by conditions 
$$b\left( xa\left( x \right) \right)=a\left( x \right),    \qquad a\left( x{{b}^{-1}}\left( x \right) \right)=b\left( x \right),$$
then}
	$$\left[ n,\to  \right]\left( 1+x{{\left( \log a\left( x \right) \right)}^{\prime }},xa\left( x \right) \right)=\left[ n,\to  \right]\left( {{b}^{n}}\left( x \right),x \right),\eqno	(1)$$
	$$\left[ n,\to  \right]\left( a\left( x \right),xa\left( x \right) \right)=\left[ n,\to  \right]\left( 1-x{{\left( \log b\left( x \right) \right)}^{\prime }}{{b}^{n+1}}\left( x \right),x \right).	\eqno(2)$$
{\bfseries Proof.} By the Lagrange inversion theorem 
$$\left[ {{x}^{n}} \right]{{a}^{m}}\left( x \right)=\frac{m}{m+n}\left[ {{x}^{n}} \right]{{b}^{m+n}}\left( x \right)=\left[ {{x}^{n}} \right]\left( 1-x{{\left( \log b\left( x \right) \right)}^{\prime }} \right){{b}^{m+n}}\left( x \right),$$
$$\left[ {{x}^{n}} \right]\left( 1+x{{\left( \log a\left( x \right) \right)}^{\prime }} \right){{a}^{m}}\left( x \right)=\frac{m+n}{m}\left[ {{x}^{n}} \right]{{a}^{m}}\left( x \right)=\left[ {{x}^{n}} \right]{{b}^{m+n}}\left( x \right).$$
Denote 
$$\left[ {{x}^{n}} \right]{{a}^{m}}\left( x \right)=a_{n}^{\left( m \right)},  \qquad\left[ {{x}^{n}} \right]\left( 1+x{{\left( \log a\left( x \right) \right)}^{\prime }} \right){{a}^{m}}\left( x \right)=c_{n}^{\left( m \right)},$$
$$\left[ {{x}^{n}} \right]{{b}^{m}}\left( x \right)=b_{n}^{\left( m \right)},  \qquad\left[ {{x}^{n}} \right]\left( 1-x{{\left( \log b\left( x \right) \right)}^{\prime }} \right){{b}^{m}}\left( x \right)=d_{n}^{\left( m \right)}.$$
Then the identities (1), (2) become obvious:
$$\left( 1+x{{\left( \log a\left( x \right) \right)}^{\prime }},xa\left( x \right) \right)=$$
$$=\left( \begin{matrix}
   c_{0}^{\left( 0 \right)} & 0 & 0 & 0 & \cdots   \\
   c_{1}^{\left( 0 \right)} & c_{0}^{\left( 1 \right)} & 0 & 0 & \cdots   \\
   c_{2}^{\left( 0 \right)} & c_{1}^{\left( 1 \right)} & c_{0}^{\left( 2 \right)} & 0 & \cdots   \\
   c_{3}^{\left( 0 \right)} & c_{2}^{\left( 1 \right)} & c_{1}^{\left( 2 \right)} & c_{0}^{\left( 3 \right)} & \cdots   \\
   \vdots  & \vdots  & \vdots  & \vdots  & \ddots   \\
\end{matrix} \right)=\left( \begin{matrix}
   b_{0}^{\left( 0 \right)} & 0 & 0 & 0 & \cdots   \\
   b_{1}^{\left( 1 \right)} & b_{0}^{\left( 1 \right)} & 0 & 0 & \cdots   \\
   b_{2}^{\left( 2 \right)} & b_{1}^{\left( 2 \right)} & b_{0}^{\left( 2 \right)} & 0 & \cdots   \\
   b_{3}^{\left( 3 \right)} & b_{2}^{\left( 3 \right)} & b_{1}^{\left( 3 \right)} & b_{0}^{\left( 3 \right)} & \cdots   \\
   \vdots  & \vdots  & \vdots  & \vdots  & \ddots   \\
\end{matrix} \right),$$
$$\left( a\left( x \right),xa\left( x \right) \right)=$$
$$=\left( \begin{matrix}
   a_{0}^{\left( 1 \right)} & 0 & 0 & 0 & \cdots   \\
   a_{1}^{\left( 1 \right)} & a_{0}^{\left( 2 \right)} & 0 & 0 & \cdots   \\
   a_{2}^{\left( 1 \right)} & a_{1}^{\left( 2 \right)} & a_{0}^{\left( 3 \right)} & 0 & \cdots   \\
   a_{3}^{\left( 1 \right)} & a_{2}^{\left( 2 \right)} & a_{1}^{\left( 3 \right)} & a_{0}^{\left( 4 \right)} & \cdots   \\
   \vdots  & \vdots  & \vdots  & \vdots  & \ddots   \\
\end{matrix} \right)=\left( \begin{matrix}
   d_{0}^{\left( 1 \right)} & 0 & 0 & 0 & \cdots   \\
   d_{1}^{\left( 2 \right)} & d_{0}^{\left( 2 \right)} & 0 & 0 & \cdots   \\
   d_{2}^{\left( 3 \right)} & d_{1}^{\left( 3 \right)} & d_{0}^{\left( 3 \right)} & 0 & \cdots   \\
   d_{3}^{\left( 4 \right)} & d_{2}^{\left( 4 \right)} & d_{1}^{\left( 4 \right)} & d_{0}^{\left( 4 \right)} & \cdots   \\
   \vdots  & \vdots  & \vdots  & \vdots  & \ddots   \\
\end{matrix} \right).$$

Note that
$${{\left( a\left( x \right),xa\left( x \right) \right)}^{-1}}=\left( {{b}^{-1}}\left( x \right),x{{b}^{-1}}\left( x \right) \right),$$
$${{\left( 1+x{{\left( \log a\left( x \right) \right)}^{\prime }},xa\left( x \right) \right)}^{-1}}=\left( 1-x{{\left( \log b\left( x \right) \right)}^{\prime }},x{{b}^{-1}}\left( x \right) \right).$$
\section{Transformations of the first type}
Let
$$b\left( xa\left( x \right) \right)=a\left( x \right),    \qquad a\left( x{{b}^{-1}}\left( x \right) \right)=b\left( x \right),$$
where
$$a\left( x \right)=\frac{1}{1-\varphi x+\beta {{x}^{2}}},  \qquad 1+x{{\left( \log a\left( x \right) \right)}^{\prime }}=\frac{1-\beta {{x}^{2}}}{1-\varphi x+\beta {{x}^{2}}},$$
$$b\left( x \right)=\frac{1+\varphi x+\sqrt{1+2\varphi x+\left( {{\varphi }^{2}}-4\beta  \right){{x}^{2}}}}{2},$$
$$\beta {{x}^{2}}{{b}^{-1}}\left( x \right)=\frac{1+\varphi x-\sqrt{1+2\varphi x+\left( {{\varphi }^{2}}-4\beta  \right){{x}^{2}}}}{2},$$
$$1-x{{\left( \log b\left( x \right) \right)}^{\prime }}=\frac{1}{\sqrt{1+2\varphi x+\left( {{\varphi }^{2}}-4\beta  \right){{x}^{2}}}}.$$
Series
$${{b}^{n}}\left( x \right)={{\left( \frac{1+\varphi x+\sqrt{1+2\varphi x+\left( {{\varphi }^{2}}-4\beta  \right){{x}^{2}}}}{2} \right)}^{n}},$$
$${{\beta }^{n}}{{x}^{2n}}{{b}^{-n}}\left( x \right)={{\left( \frac{1+\varphi x-\sqrt{1+2\varphi x+\left( {{\varphi }^{2}}-4\beta  \right){{x}^{2}}}}{2} \right)}^{n}}$$
can be represented as
$${{b}^{n}}\left( x \right)=\frac{{{c}_{n}}\left( \varphi ,\beta ,x \right)+{{s}_{n}}\left( \varphi ,\beta ,x \right)\sqrt{1+2\varphi x+\left( {{\varphi }^{2}}-4\beta  \right){{x}^{2}}}}{2},$$
$${{\beta }^{n}}{{x}^{2n}}{{b}^{-n}}\left( x \right)=\frac{{{c}_{n}}\left( \varphi ,\beta ,x \right)-{{s}_{n}}\left( \varphi ,\beta ,x \right)\sqrt{1+2\varphi x+\left( {{\varphi }^{2}}-4\beta  \right){{x}^{2}}}}{2},$$
where ${{c}_{n}}\left( \varphi ,\beta ,x \right)$ is the polynomial of degree $\le n$, ${{s}_{n}}\left( \varphi ,\beta ,x \right)$ is the polynomial of degree  $<n$, and
$${{s}_{n}}\left( \varphi ,\beta ,x \right)\sqrt{1-2\varphi +\left( {{\varphi }^{2}}-4\beta  \right){{x}^{2}}}=\sqrt{c_{n}^{2}\left( \varphi ,\beta ,x \right)-4{{\beta }^{n}}{{x}^{2n}}},$$
as it follows from
$${{b}^{2n}}\left( x \right)-{{c}_{n}}\left( \varphi ,\beta ,x \right){{b}^{n}}\left( x \right)+{{\beta }^{n}}{{x}^{2n}}=0.$$
Let ${{J}_{n}}$ is the operator exchanging the coefficients of the polynomial of degree $n$ in reverse order: ${{J}_{n}}c\left( x \right)={{x}^{n}}c\left( {1}/{x}\; \right)$, where $c\left( x \right)$ is the polynomial of degree $\le n$. Comparison the identities 
$${{b}^{n}}\left( x \right)={{c}_{n}}\left( \varphi ,\beta ,x \right)-{{\beta }^{n}}{{x}^{2n}}{{b}^{-n}}\left( x \right),$$
$$\left( 1-x{{\left( \log b\left( x \right) \right)}^{\prime }} \right){{b}^{n}}\left( x \right)={{s}_{n}}\left( \varphi ,\beta ,x \right)+\frac{{{\beta }^{n}}{{x}^{2n}}{{b}^{-n}}\left( x \right)}{\sqrt{1+2\varphi x+\left( {{\varphi }^{2}}-4\beta  \right){{x}^{2}}}},$$
respectively with the identities (1), (2) gives:
$${{c}_{0}}\left( \varphi ,\beta ,x \right)=2,   \qquad{{c}_{n}}\left( \varphi ,\beta ,x \right)={{J}_{n}}{{D}_{n}}\left( x+\varphi ,\beta  \right),$$   
$${{s}_{0}}\left( \varphi ,\beta ,x \right)=0,   \qquad{{s}_{n}}\left( \varphi ,\beta ,x \right)={{J}_{n-1}}{{E}_{n-1}}\left( x+\varphi ,\beta  \right).$$
{\bfseries Example 1.} If $\varphi =0$, $\beta =1$, then
$$b\left( x \right)=\frac{1+\sqrt{1-4{{x}^{2}}}}{2},  \qquad{{b}^{-1}}\left( x \right)=\frac{1-\sqrt{1-4{{x}^{2}}}}{2{{x}^{2}}}=C\left( {{x}^{2}} \right),$$
$${{C}^{n}}\left( {{x}^{2}} \right)=\frac{{{x}^{n}}{{C}_{n}}\left( {1}/{x}\; \right)-{{x}^{n-1}}{{S}_{n-1}}\left( {1}/{x}\; \right)\sqrt{1-4{{x}^{2}}}}{2{{x}^{2n}}}=$$
$$=\frac{-{{x}^{n-2}}{{S}_{n-2}}\left( {1}/{x}\; \right)+{{x}^{n-1}}{{S}_{n-1}}\left( {1}/{x}\; \right)C\left( {{x}^{2}} \right)}{{{x}^{2n-2}}},$$
$${{C}^{-n}}\left( {{x}^{2}} \right)=\frac{{{x}^{n}}{{C}_{n}}\left( {1}/{x}\; \right)+{{x}^{n-1}}{{S}_{n-1}}\left( {1}/{x}\; \right)\sqrt{1-4{{x}^{2}}}}{2}=$$
$$={{x}^{n}}{{S}_{n}}\left( {1}/{x}\; \right)-{{x}^{n+1}}{{S}_{n-1}}\left( {1}/{x}\; \right)C\left( {{x}^{2}} \right),$$
where we used identities
$${{C}_{n}}\left( x \right)=x{{S}_{n-1}}\left( x \right)-2{{S}_{n-2}}\left( x \right)=2{{S}_{n}}\left( x \right)-x{{S}_{n-1}}\left( x \right),$$ 
or
$${{x}^{n}}{{C}_{n}}\left( {1}/{x}\; \right)={{x}^{n-1}}{{S}_{n-1}}\left( {1}/{x}\; \right)-2{{x}^{n}}{{S}_{n-2}}\left( {1}/{x}\; \right)=2{{x}^{n}}{{S}_{n}}\left( {1}/{x}\; \right)-{{x}^{n-1}}{{S}_{n-1}}\left( {1}/{x}\; \right).$$
Denote ${{S}_{-n}}\left( x \right)=-{{S}_{n-2}}\left( x \right)$. Then
$${{C}^{k}}\left( {{x}^{2}} \right)=-{{\left( \frac{1}{x} \right)}^{k}}{{S}_{k-2}}\left( \frac{1}{x} \right)+{{\left( \frac{1}{x} \right)}^{k-1}}{{S}_{k-1}}\left( \frac{1}{x} \right)C\left( {{x}^{2}} \right),$$
$k=0$, $\pm 1$, $\pm 2$, … . Hence, as shown in [16],
$${{C}^{k}}\left( x \right)=-{{\left( \frac{1}{\sqrt{x}} \right)}^{k}}{{S}_{k-2}}\left( \frac{1}{\sqrt{x}} \right)+{{\left( \frac{1}{\sqrt{x}} \right)}^{k-1}}{{S}_{k-1}}\left( \frac{1}{\sqrt{x}} \right)C\left( x \right).$$ 

Note that
$$\left( \frac{1-\beta {{x}^{2}}}{1-\varphi x+\beta {{x}^{2}}},\frac{x}{1-\varphi x+\beta {{x}^{2}}} \right){{P}^{-2\varphi }}=M\left( \frac{1-\beta {{x}^{2}}}{1-\varphi x+\beta {{x}^{2}}},\frac{x}{1-\varphi x+\beta {{x}^{2}}} \right)M,$$
$$\left( \frac{1}{1-\varphi x+\beta {{x}^{2}}},\frac{x}{1-\varphi x+\beta {{x}^{2}}} \right){{P}^{-2\varphi }}=M\left( \frac{1}{1-\varphi x+\beta {{x}^{2}}},\frac{x}{1-\varphi x+\beta {{x}^{2}}} \right)M.$$
Respectively,
$${{P}^{2\varphi }}\left( 1-x{{\left( \log b\left( x \right) \right)}^{\prime }},x{{b}^{-1}}\left( x \right) \right)=M\left( 1-x{{\left( \log b\left( x \right) \right)}^{\prime }},x{{b}^{-1}}\left( x \right) \right)M,$$
$${{P}^{2\varphi }}\left( {{b}^{-1}}\left( x \right),x{{b}^{-1}}\left( x \right) \right)=M\left( {{b}^{-1}}\left( x \right),x{{b}^{-1}}\left( x \right) \right)M.$$
\section{Transformations of the second type}
Denote
$${{x}^{n}}-1=\prod\limits_{m=0}^{n-1}{\left( x-{{e}^{\left( n,\text{ }m \right)}} \right)},  \qquad{{x}^{n}}+1=\frac{{{x}^{2n}}-1}{{{x}^{n}}-1}=\prod\limits_{m=1}^{n}{\left( x-{{e}^{\left( 2n,\text{ }2m-1 \right)}} \right)},$$
$${{e}^{\left( n,\text{ }m \right)}}=\cos \frac{2\pi m}{n}+i\sin \frac{2\pi m}{n},  \qquad i=\sqrt{-1}.$$
Then
$$\frac{{{\left( x+1 \right)}^{n}}+{{\left( x-1 \right)}^{n}}}{2}={{2}^{n-1}}\left( 1,{x}/{2}\; \right){{\left( {{P}^{-{1}/{2}\;}} \right)}^{T}}{{J}_{n}}{{P}^{T}}\left( {{x}^{n}}+1 \right)=$$
$$=\prod\limits_{m=1}^{n}{\left( x+\frac{1+{{e}^{\left( 2n,\text{ }2m-1 \right)}}}{1-{{e}^{\left( 2n,\text{ }2m-1 \right)}}} \right)}=\prod\limits_{m=1}^{n}{\left( x+i\text{ctg}\frac{2m-1}{2n}\pi  \right)},$$
$$\frac{{{\left( x+1 \right)}^{n}}-{{\left( x-1 \right)}^{n}}}{2}={{2}^{n-1}}\left( 1,{x}/{2}\; \right){{\left( {{P}^{-{1}/{2}\;}} \right)}^{T}}{{J}_{n}}{{P}^{T}}\left( {{x}^{n}}-1 \right)=$$
$$=n\prod\limits_{m=1}^{n-1}{\left( x+\frac{1+{{e}^{\left( n,\text{ }m \right)}}}{1-{{e}^{\left( n,\text{ }m \right)}}} \right)}=n\prod\limits_{m=1}^{n-1}{\left( x+i\text{ctg}\frac{m}{n}\pi  \right)},$$
Consider matrices composed of the even and odd  columns of the Pascal matrix:
$$\left( \frac{1-x}{{{\left( 1-x \right)}^{2}}},\frac{{{x}^{2}}}{{{\left( 1-x \right)}^{2}}} \right)=\left( \begin{matrix}
   1 & 0 & 0 & \cdots   \\
   1 & 0 & 0 & \cdots   \\
   1 & 1 & 0 & \cdots   \\
   1 & 3 & 0 & \cdots   \\
   1 & 6 & 1 & \cdots   \\
   1 & 10 & 5 & \cdots   \\
   \vdots  & \vdots  & \vdots  & \ddots   \\
\end{matrix} \right), \quad\left( \frac{x}{{{\left( 1-x \right)}^{2}}},\frac{{{x}^{2}}}{{{\left( 1-x \right)}^{2}}} \right)=\left( \begin{matrix}
   0 & 0 & 0 & \cdots   \\
   1 & 0 & 0 & \cdots   \\
   2 & 0 & 0 & \cdots   \\
   3 & 1 & 0 & \cdots   \\
   4 & 4 & 0 & \cdots   \\
   5 & 10 & 1 & \cdots   \\
   \vdots  & \vdots  & \vdots  & \ddots   \\
\end{matrix} \right).$$
Denote
$$\left[ n,\to  \right]\left( \frac{1-x}{{{\left( 1-x \right)}^{2}}},\frac{{{x}^{2}}}{{{\left( 1-x \right)}^{2}}} \right)={{\tilde{t}}_{n}}\left( 1,1,x \right),
\quad\left[ n,\to  \right]\left( \frac{x}{{{\left( 1-x \right)}^{2}}},\frac{{{x}^{2}}}{{{\left( 1-x \right)}^{2}}} \right)={{\tilde{u}}_{n}}\left( 1,1,x \right).$$
Since
$${{\tilde{t}}_{n}}\left( 1,1,{{x}^{2}} \right)=\frac{{{\left( 1+x \right)}^{n}}+{{\left( 1-x \right)}^{n}}}{2}={{J}_{n}}\frac{{{\left( x+1 \right)}^{n}}+{{\left( x-1 \right)}^{n}}}{2}=$$
$$={{d}_{n}}\prod\limits_{m=1}^{\left\lfloor {n}/{2}\; \right\rfloor }{\left( {{x}^{2}}+\text{t}{{\text{g}}^{2}}\frac{2m-1}{2n}\pi  \right)},  \qquad{{d}_{n}}=1,\quad n,$$ 
(first value of ${{d}_{n}}$ is taken for even $n$, the second taken for odd $n$; we will agree to use this rule for other similar quantities),
$${{\tilde{u}}_{n}}\left( 1,1,{{x}^{2}} \right)=\frac{{{\left( 1+x \right)}^{n}}-{{\left( 1-x \right)}^{n}}}{2x}={{J}_{n}}\frac{{{\left( x+1 \right)}^{n}}-{{\left( x-1 \right)}^{n}}}{2}=$$
$$={{l}_{n}}\prod\limits_{m=1}^{\left\lfloor {\left( n-1 \right)}/{2}\; \right\rfloor }{\left( {{x}^{2}}+\text{t}{{\text{g}}^{2}}\frac{m}{n}\pi  \right)}, \qquad{{l}_{n}}=n,\quad 1,$$
then
$${{\tilde{t}}_{n}}\left( 1,1,x \right)={{d}_{n}}\prod\limits_{m=1}^{\left\lfloor {n}/{2}\; \right\rfloor }{\left( x+\text{t}{{\text{g}}^{2}}\frac{2m-1}{2n}\pi  \right)}, \quad{{\tilde{u}}_{n}}\left( 1,1,x \right)={{l}_{n}}\prod\limits_{m=1}^{\left\lfloor {\left( n-1 \right)}/{2}\; \right\rfloor }{\left( x+\text{t}{{\text{g}}^{2}}\frac{m}{n}\pi  \right)}.$$
Then
$$\left( \frac{1-x}{{{\left( 1-x \right)}^{2}}},\frac{{{x}^{2}}}{{{\left( 1-x \right)}^{2}}} \right)\left( \frac{1}{1-x},\frac{x}{1-x} \right)=\left( \frac{1-x}{1-2x},\frac{{{x}^{2}}}{1-2x} \right),$$
$$\left( \frac{x}{{{\left( 1-x \right)}^{2}}},\frac{{{x}^{2}}}{{{\left( 1-x \right)}^{2}}} \right)\left( \frac{1}{1-x},\frac{x}{1-x} \right)=\left( \frac{x}{1-2x},\frac{{{x}^{2}}}{1-2x} \right),$$
$$\left( \frac{1-x}{1-2x},\frac{{{x}^{2}}}{1-2x} \right)=\left( \begin{matrix}
   1 & 0 & 0 & \cdots   \\
   1 & 0 & 0 & \cdots   \\
   2 & 1 & 0 & \cdots   \\
   4 & 3 & 0 & \cdots   \\
   8 & 8 & 1 & \cdots   \\
   16 & 20 & 5 & \cdots   \\
   \vdots  & \vdots  & \vdots  & \ddots   \\
\end{matrix} \right), \quad\left( \frac{x}{1-2x},\frac{{{x}^{2}}}{1-2x} \right)=\left( \begin{matrix}
   0 & 0 & 0 & \cdots   \\
   1 & 0 & 0 & \cdots   \\
   2 & 0 & 0 & \cdots   \\
   4 & 1 & 0 & \cdots   \\
   8 & 4 & 0 & \cdots   \\
   16 & 12 & 1 & \cdots   \\
   \vdots  & \vdots  & \vdots  & \ddots   \\
\end{matrix} \right).$$
Denote
$$\left[ n,\to  \right]\left( \frac{1-x}{1-2x},\frac{{{x}^{2}}}{1-2x} \right)={{\tilde{t}}_{n}}\left( 1,0,x \right),   \quad\left[ n,\to  \right]\left( \frac{x}{1-2x},\frac{{{x}^{2}}}{1-2x} \right)={{\tilde{u}}_{n}}\left( 1,0,x \right).$$
Then
$${{\tilde{t}}_{n}}\left( 1,0,x \right)={{d}_{n}}\prod\limits_{m=1}^{\left\lfloor {n}/{2}\; \right\rfloor }{\left( x+{{\sec }^{2}}\frac{2m-1}{2n}\pi  \right)},
\quad{{\tilde{u}}_{n}}\left( 1,0,x \right)={{l}_{n}}\prod\limits_{m=1}^{\left\lfloor {\left( n-1 \right)}/{2}\; \right\rfloor }{\left( x+{{\sec }^{2}}\frac{m}{n}\pi  \right)}.$$
Hence, $n$th row of the matrix 
$$\left( \frac{1-\varphi x}{1-2\varphi x+\beta {{x}^{2}}},\frac{{{x}^{2}}}{1-2\varphi x+\beta {{x}^{2}}} \right)=\left( \frac{1-\varphi x}{1-2\varphi x},\frac{{{x}^{2}}}{1-2\varphi x} \right)\left( \frac{1}{1+\beta x},\frac{x}{1+\beta x} \right)$$
is the polynomial
$${{\tilde{t}}_{n}}\left( \varphi ,\beta ,x \right)={{p}_{n}}\prod\limits_{m=1}^{\left\lfloor {n}/{2}\; \right\rfloor }{\left( x+\frac{{{\varphi }^{2}}-\beta {{\cos }^{2}}\frac{2m-1}{2n}\pi }{{{\cos }^{2}}\frac{2m-1}{2n}\pi } \right)}, \qquad{{p}_{n}}=1,\quad n\varphi ;$$
$n$th row of the matrix 
$$\left( \frac{x}{1-2\varphi x+\beta {{x}^{2}}},\frac{{{x}^{2}}}{1-2\varphi x+\beta {{x}^{2}}} \right)=\left( \frac{x}{1-2\varphi x},\frac{{{x}^{2}}}{1-2\varphi x} \right)\left( \frac{1}{1+\beta x},\frac{x}{1+\beta x} \right)$$
is the polynomial
$${{\tilde{u}}_{n}}\left( \varphi ,\beta ,x \right)={{r}_{n}}\prod\limits_{m=1}^{\left\lfloor {\left( n-1 \right)}/{2}\; \right\rfloor }{\left( x+\frac{{{\varphi }^{2}}-\beta {{\cos }^{2}}\frac{m}{n}\pi }{{{\cos }^{2}}\frac{m}{n}\pi } \right)}, \qquad{{r}_{n}}=n\varphi , \quad1.$$

Sequences of columns of the matrices 
$$\left( \frac{1-\varphi x}{1-2\varphi x+\beta {{x}^{2}}},\frac{{{x}^{2}}}{1-2\varphi x+\beta {{x}^{2}}} \right), \quad\left( \frac{x}{1-2\varphi x+\beta {{x}^{2}}},\frac{{{x}^{2}}}{1-2\varphi x+\beta {{x}^{2}}} \right)$$
coincides respectively with the sequences of even and odd columns of the matrices
$$\left( \frac{1-\varphi x}{1-2\varphi x+\beta {{x}^{2}}},\frac{x}{\sqrt{\left( 1-2\varphi x+\beta {{x}^{2}} \right)}} \right), \quad\left( \frac{1}{\sqrt{1-2\varphi x+\beta {{x}^{2}}}},\frac{x}{\sqrt{1-2\varphi x+\beta {{x}^{2}}}} \right).$$
Let
$$b\left( xa\left( x \right) \right)=a\left( x \right),    \qquad a\left( x{{b}^{-1}}\left( x \right) \right)=b\left( x \right),$$
where
$$a\left( x \right)=\frac{1}{\sqrt{1-2\varphi x+\beta {{x}^{2}}}}, \qquad1+x{{\left( \log a\left( x \right) \right)}^{\prime }}=\frac{1-\varphi x}{1-2\varphi x+\beta {{x}^{2}}},$$
$$b\left( x \right)=\varphi x+\sqrt{1+\left( {{\varphi }^{2}}-\beta  \right){{x}^{2}}},  \qquad{{b}^{-1}}\left( x \right)=\frac{\varphi x-\sqrt{1+\left( {{\varphi }^{2}}-\beta  \right){{x}^{2}}}}{\beta {{x}^{2}}-1},$$
$$1-x{{\left( \log b\left( x \right) \right)}^{\prime }}=\frac{1}{b\left( x \right)\sqrt{1+\left( {{\varphi }^{2}}-\beta  \right){{x}^{2}}}}.$$
Series
$${{b}^{n}}\left( x \right)={{\left( \varphi x+\sqrt{1+\left( {{\varphi }^{2}}-\beta  \right){{x}^{2}}} \right)}^{n}},
\quad{{\left( \beta {{x}^{2}}-1 \right)}^{n}}{{b}^{-n}}\left( x \right)={{\left( \varphi x-\sqrt{1+\left( {{\varphi }^{2}}-\beta  \right){{x}^{2}}} \right)}^{n}},$$ 
can be represented as
$${{b}^{n}}\left( x \right)={{t}_{n}}\left( \varphi ,\beta ,x \right)+{{u}_{n}}\left( \varphi ,\beta ,x \right)\sqrt{1+\left( {{\varphi }^{2}}-\beta  \right){{x}^{2}}},$$
$${{\left( \beta {{x}^{2}}-1 \right)}^{n}}{{b}^{-n}}\left( x \right)={{t}_{n}}\left( \varphi ,\beta ,x \right)-{{u}_{n}}\left( \varphi ,\beta ,x \right)\sqrt{1+\left( {{\varphi }^{2}}-\beta  \right){{x}^{2}}},$$
where ${{t}_{n}}\left( \varphi ,\beta ,x \right)$ is the polynomial of degree $\le n$, 
$$\left[ {{x}^{2n}} \right]{{t}_{2m+1}}\left( \varphi ,\beta ,x \right)=0,   \qquad\left[ {{x}^{2n+1}} \right]{{t}_{2m}}\left( \varphi ,\beta ,x \right)=0;$$
${{u}_{n}}\left( \varphi ,\beta ,x \right)$ is the polynomial of degree $<n$, 
$$\left[ {{x}^{2n}} \right]{{u}_{2m}}\left( \varphi ,\beta ,x \right)=0,  \qquad\left[ {{x}^{2n+1}} \right]{{u}_{2m+1}}\left( \varphi ,\beta ,x \right)=0,$$
and
$${{u}_{n}}\left( \varphi ,\beta ,x \right)\sqrt{1+\left( {{\varphi }^{2}}-\beta  \right){{x}^{2}}}=\sqrt{t_{n}^{2}\left( \varphi ,\beta ,x \right)-{{\left( \beta {{x}^{2}}-1 \right)}^{n}}},$$
as it follows from
$${{b}^{2n}}\left( x \right)-2{{t}_{n}}\left( \varphi ,\beta ,x \right){{b}^{n}}\left( x \right)+{{\left( \beta {{x}^{2}}-1 \right)}^{n}}=0.$$
Comparison the identities 
$${{b}^{n}}\left( x \right)={{t}_{n}}\left( \varphi ,\beta ,x \right)+{{u}_{n}}\left( \varphi ,\beta ,x \right)\sqrt{1+\left( {{\varphi }^{2}}-\beta  \right){{x}^{2}}},$$
$$\left( 1-x{{\left( \log b\left( x \right) \right)}^{\prime }} \right){{b}^{n+1}}\left( x \right)={{u}_{n}}\left( \varphi ,\beta ,x \right)+\frac{{{t}_{n}}\left( \varphi ,\beta ,x \right)}{\sqrt{1+\left( {{\varphi }^{2}}-\beta  \right){{x}^{2}}}}$$
respectively with the identities (1), (2) gives:
$${{t}_{n}}\left( \varphi ,\beta ,x \right)={{J}_{n}}{{\tilde{t}}_{n}}\left( \varphi ,\beta ,{{x}^{2}} \right)=$$ 
$$={{p}_{n}}\prod\limits_{m=1}^{\left\lfloor {n}/{2}\; \right\rfloor }{\left( \frac{{{\varphi }^{2}}-\beta {{\cos }^{2}}\frac{2m-1}{2n}\pi }{{{\cos }^{2}}\frac{2m-1}{2n}\pi } \right)}\prod\limits_{m=1}^{n}{\left( x+\frac{i\cos \frac{2m-1}{2n}\pi }{\sqrt{{{\varphi }^{2}}-\beta {{\cos }^{2}}\frac{2m-1}{2n}\pi }} \right)},$$
where if ${{\varphi }^{2}}-\beta {{\cos }^{2}}\frac{2m-1}{2n}\pi =0$, we accept
$$\frac{{{\varphi }^{2}}-\beta {{\cos }^{2}}\frac{2m-1}{2n}\pi }{{{\cos }^{2}}\frac{2m-1}{2n}\pi }=1,  \qquad x+\frac{i\cos \frac{2m-1}{2n}\pi }{\sqrt{{{\varphi }^{2}}-\beta {{\cos }^{2}}\frac{2m-1}{2n}\pi }}=1;$$
$${{u}_{n}}\left( \varphi ,\beta ,x \right)={{J}_{n-1}}{{\tilde{u}}_{n}}\left( \varphi ,\beta ,{{x}^{2}} \right)=$$
$$={{r}_{n}}\prod\limits_{m=1}^{\left\lfloor {\left( n-1 \right)}/{2}\; \right\rfloor }{\left( \frac{{{\varphi }^{2}}-\beta {{\cos }^{2}}\frac{m}{n}\pi }{{{\cos }^{2}}\frac{m}{n}\pi } \right)}\prod\limits_{m=1}^{n-1}{\left( x+\frac{i\cos \frac{m}{n}\pi }{\sqrt{{{\varphi }^{2}}-\beta {{\cos }^{2}}\frac{m}{n}\pi }} \right)},$$
where if  ${{\varphi }^{2}}-\beta {{\cos }^{2}}\frac{m}{n}\pi =0$, we accept
$$\frac{{{\varphi }^{2}}-\beta {{\cos }^{2}}\frac{m}{n}\pi }{{{\cos }^{2}}\frac{m}{n}\pi }=1,  \qquad x+\frac{i\cos \frac{m}{n}\pi }{\sqrt{{{\varphi }^{2}}-\beta {{\cos }^{2}}\frac{m}{n}\pi }}=1.$$ 

We will take into account the identities 
$$\prod\limits_{m=1}^{\left\lfloor {n}/{2}\; \right\rfloor }{{{\cos }^{2}}\frac{2m-1}{2n}\pi }=\frac{1}{{{2}^{n-1}}}, \quad\frac{n}{{{2}^{n-1}}}; \qquad\prod\limits_{m=1}^{\left\lfloor {\left( n-1 \right)}/{2}\; \right\rfloor }{{{\cos }^{2}}\frac{m}{n}\pi }=\frac{n}{{{2}^{n-1}}},\quad  \frac{1}{{{2}^{n-1}}};$$
$$\prod\limits_{m=1}^{\left\lfloor {n}/{2}\; \right\rfloor }{{{\sin }^{2}}\frac{2m-1}{2n}\pi }=\frac{1}{{{2}^{n-1}}},   \qquad\prod\limits_{m=1}^{\left\lfloor {\left( n-1 \right)}/{2}\; \right\rfloor }{{{\sin }^{2}}\frac{m}{n}\pi }=\frac{n}{{{2}^{n-1}}},$$
and note the cases:
$${{t}_{n}}\left( 1,1,x \right)=\frac{{{\left( x+1 \right)}^{n}}+{{\left( x-1 \right)}^{n}}}{2},   \qquad{{u}_{n}}\left( 1,1,x \right)=\frac{{{\left( x+1 \right)}^{n}}-{{\left( x-1 \right)}^{n}}}{2};$$
$${{t}_{n}}\left( 1,0,x \right)=\frac{{{\left( x+\sqrt{1+{{x}^{2}}} \right)}^{n}}+{{\left( x-\sqrt{1+{{x}^{2}}} \right)}^{n}}}{2},$$
$${{u}_{n}}\left( 1,0,x \right)=\frac{{{\left( x+\sqrt{1+{{x}^{2}}} \right)}^{n}}-{{\left( x-\sqrt{1+{{x}^{2}}} \right)}^{n}}}{2\sqrt{1+{{x}^{2}}}};$$
$${{t}_{2n}}\left( 0.-1,x \right)={{\left( 1+{{x}^{2}} \right)}^{n}},    \qquad{{u}_{2n}}\left( 0,-1,x \right)=0,$$
$${{t}_{2n+1}}\left( 0,-1,x \right)=0,     \qquad{{u}_{2n+1}}\left( 0,-1,x \right)={{\left( 1+{{x}^{2}} \right)}^{n}};$$
$${{t}_{2n}}\left( 0,0,x \right)=1,       \qquad{{u}_{2n}}\left( 0,0,x \right)=0,$$
$${{t}_{2n+1}}\left( 0,0,x \right)=0,      \qquad{{u}_{2n+1}}\left( 0,0,x \right)=1.$$

Note that
$${{P}^{-2\varphi }}\left( \frac{1-\varphi x}{1-2\varphi x+\beta {{x}^{2}}},\frac{x}{\sqrt{1-2\varphi x+\beta {{x}^{2}}}} \right)
=M\left( \frac{1-\varphi x}{1-2\varphi x+\beta {{x}^{2}}},\frac{x}{\sqrt{1-2\varphi x+\beta {{x}^{2}}}} \right)M,$$
$${{P}^{-2\varphi }}\left( \frac{1}{\sqrt{1-2\varphi x+\beta {{x}^{2}}}},\frac{x}{\sqrt{1-2\varphi x+\beta {{x}^{2}}}} \right)
=M\left( \frac{1}{\sqrt{1-2\varphi x+\beta {{x}^{2}}}},\frac{x}{\sqrt{1-2\varphi x+\beta {{x}^{2}}}} \right)M.$$
Respectively,
$$\left( 1-x{{\left( \log b\left( x \right) \right)}^{\prime }},x{{b}^{-1}}\left( x \right) \right){{P}^{2\varphi }}=M\left( 1-x{{\left( \log b\left( x \right) \right)}^{\prime }},x{{b}^{-1}}\left( x \right) \right)M,$$
$$\left( {{b}^{-1}}\left( x \right),x{{b}^{-1}}\left( x \right) \right){{P}^{2\varphi }}=M\left( {{b}^{-1}}\left( x \right),x{{b}^{-1}}\left( x \right) \right)M.$$
\section{Fibonacci bases}
Let
$$b\left( xa\left( x \right) \right)=a\left( x \right),    \qquad a\left( x{{b}^{-1}}\left( x \right) \right)=b\left( x \right),$$
where
$$a\left( x \right)=\frac{x}{2}+\sqrt{1+\frac{{{x}^{2}}}{4}},$$
$$ b\left( x \right)=\sqrt{1+x},    \qquad1-x{{\left( \log b\left( x \right) \right)}^{\prime }}=\left( 1+\frac{x}{2} \right)\frac{1}{1+x}.$$
Comparison the identities 
$${{b}^{2n}}\left( x \right)={{\left( 1+x \right)}^{n}} ,  \qquad\left( 1-x{{\left( \log b\left( x \right) \right)}^{\prime }} \right){{b}^{2n+2}}\left( x \right)=\left( 1+\frac{x}{2} \right){{\left( 1+x \right)}^{n}}$$
respectively with the identities (1), (2) gives:
$$\left[ 2n,\to  \right]\left( 1+x{{\left( \log a\left( x \right) \right)}^{\prime }},xa\left( x \right) \right)={{x}^{n}}{{\left( 1+x \right)}^{n}}=\left( 1,x\left( 1+x \right) \right){{x}^{n}},$$
$$\left[ 2n+1,\to  \right]\left( a\left( x \right),xa\left( x \right) \right)=\left( \frac{1}{2}+x \right){{x}^{n}}{{\left( 1+x \right)}^{n}}=\frac{1}{2}\left( 1+2x,x\left( 1+x \right) \right){{x}^{n}},$$
$$\left( 1,x\left( 1+x \right) \right)=\left( \begin{matrix}
   1 & 0 & 0 & 0 & 0 & 0 & \cdots   \\
   0 & 1 & 0 & 0 & 0 & 0 & \cdots   \\
   0 & 1 & 1 & 0 & 0 & 0 & \cdots   \\
   0 & 0 & 2 & 1 & 0 & 0 & \cdots   \\
   0 & 0 & 1 & 3 & 1 & 0 & \cdots   \\
   0 & 0 & 0 & 3 & 4 & 1 & \cdots   \\
   \vdots  & \vdots  & \vdots  & \vdots  & \vdots  & \vdots  & \ddots   \\
\end{matrix} \right), \\\left( 1+2x,x\left( 1+x \right) \right)=\left( \begin{matrix}
   1 & 0 & 0 & 0 & 0 & 0 & \cdots   \\
   2 & 1 & 0 & 0 & 0 & 0 & \cdots   \\
   0 & 3 & 1 & 0 & 0 & 0 & \cdots   \\
   0 & 2 & 4 & 1 & 0 & 0 & \cdots   \\
   0 & 0 & 5 & 5 & 1 & 0 & \cdots   \\
   0 & 0 & 2 & 9 & 6 & 1 & \cdots   \\
   \vdots  & \vdots  & \vdots  & \vdots  & \vdots  & \vdots  & \ddots   \\
\end{matrix} \right).$$

Let $a\left( x \right)\to \left( {{a}_{0}},{{a}_{1}},{{a}_{2}},... \right)$ mean that the series $a\left( x \right)$ is the generating function of the sequence $\left( {{a}_{0}},{{a}_{1}},{{a}_{2}},... \right)$. Note that
$$\left( 1,x\left( 1+x \right) \right)\frac{1}{1-x}=\frac{1}{1-x-{{x}^{2}}}\to \left( 1,\text{ }1,\text{ }2,\text{ }3,\text{ }5,\text{ }... \right),$$
$$\left( 1+2x,x\left( 1+x \right) \right)\frac{1}{1-x}=\frac{1+2x}{1-x-{{x}^{2}}}\to \left( 1,\text{ }3,\text{ }4,\text{ }7,\text{ }11,\text{ }... \right).$$
In this connection, matrices $\left( 1,x\left( 1+x \right) \right)$, $\left( 1+2x,x\left( 1+x \right) \right)$ are called Fibonacci and Lucas  matrices. Note that
$${{\left( {{P}^{-1}} \right)}^{T}}\left( 1,x\left( 1+x \right) \right)=\left( 1,x\left( x-1 \right) \right)=M\left( 1,x\left( 1+x \right) \right),$$
$${{\left( {{P}^{-1}} \right)}^{T}}\left( 1+2x,x\left( 1+x \right) \right)=\left( 2x-1,x\left( x-1 \right) \right)=-M\left( 1+2x,x\left( 1+x \right) \right).$$

We will build the pseudo-eigenbases of the transformations ${{\left( {{P}^{-1}} \right)}^{T}}$, the even columns of which are the columns of the matrix $\left( 1,x\left( 1+x \right) \right)$ and the odd columns are the columns of the matrix $\left( 1+2x,x\left( 1+x \right) \right)$:
$$B=\text{ }\left( \begin{matrix}
   1 & 1 & 0 & 0 & 0 & 0 & 0 & \cdots   \\
   0 & 2 & 1 & 1 & 0 & 0 & 0 & \cdots   \\
   0 & 0 & 1 & 3 & 1 & 1 & 0 & \cdots   \\
   0 & 0 & 0 & 2 & 2 & 4 & 1 & \cdots   \\
   0 & 0 & 0 & 0 & 1 & 5 & 3 & \cdots   \\
   0 & 0 & 0 & 0 & 0 & 2 & 3 & \cdots   \\
   0 & 0 & 0 & 0 & 0 & 0 & 1 & \cdots   \\
   \vdots  & \vdots  & \vdots  & \vdots  & \vdots  & \vdots  & \vdots  & \ddots   \\
\end{matrix} \right).$$

Consider matrix consisting of even columns of the matrix 
$$2{{\left( 1+x{{\left( \log a\left( x \right) \right)}^{\prime }},xa\left( x \right) \right)}^{-1}}=2\left( \left( 1+\frac{x}{2} \right)\frac{1}{1+x},\frac{1}{\sqrt{1+x}} \right)$$
and matrix consisting of odd columns of the matrix
$${{\left( a\left( x \right),xa\left( x \right) \right)}^{-1}}=\left( \frac{1}{\sqrt{1+x}},\frac{x}{\sqrt{1+x}} \right):$$
$$\left( \frac{2+x}{1+x},\frac{{{x}^{2}}}{1+x} \right)=\left( \begin{matrix}
   2 & 0 & 0 & 0 & \cdots   \\
   -1 & 0 & 0 & 0 & \cdots   \\
   1 & 2 & 0 & 0 & \cdots   \\
   -1 & -3 & 0 & 0 & \cdots   \\
   1 & 4 & 2 & 0 & \cdots   \\
   -1 & -5 & -5 & 0 & \cdots   \\
   1 & 6 & 9 & 2 & \cdots   \\
   -1 & -7 & -14 & -7 & \cdots   \\
   \vdots  & \vdots  & \vdots  & \vdots  & \ddots   \\
\end{matrix} \right), \\\left( \frac{x}{1+x},\frac{{{x}^{2}}}{1+x} \right)=\left( \begin{matrix}
   0 & 0 & 0 & 0 & \cdots   \\
   1 & 0 & 0 & 0 & \cdots   \\
   -1 & 0 & 0 & 0 & \cdots   \\
   1 & 1 & 0 & 0 & \cdots   \\
   -1 & -2 & 0 & 0 & \cdots   \\
   1 & 3 & 1 & 0 & \cdots   \\
   -1 & -4 & -3 & 0 & \cdots   \\
   1 & 5 & 6 & 1 & \cdots   \\
   \vdots  & \vdots  & \vdots  & \vdots  & \ddots   \\
\end{matrix} \right).$$
Note that
 $$\left( \frac{2+x}{1+x},\frac{{{x}^{2}}}{1+x} \right)\frac{1}{1-x}=\frac{2+x}{1+x-{{x}^{2}}}\to \left( 2,-1,\text{ }3,-4,\text{ }7,\text{ }... \right),$$
$$\left( \frac{x}{1+x},\frac{{{x}^{2}}}{1+x} \right)\frac{1}{1-x}\to \frac{x}{1+x-{{x}^{2}}}=\left( 0,\text{ }1,-1,\text{ }2,-3,\text{ }5,\text{ }... \right);$$ 
$$P\left( \frac{2+x}{1+x},\frac{{{x}^{2}}}{1+x} \right)=\left( \frac{2-x}{1-x},\frac{{{x}^{2}}}{1-x} \right)=M\left( \frac{2+x}{1+x},\frac{{{x}^{2}}}{1+x} \right),$$
$$P\left( \frac{x}{1+x},\frac{{{x}^{2}}}{1+x} \right)=\left( \frac{x}{1-x},\frac{{{x}^{2}}}{1-x} \right)=-M\left( \frac{x}{1+x},\frac{{{x}^{2}}}{1+x} \right).$$

We will build the pseudo-eigenbases of the transformations $P$, the even columns of which are the columns of the matrix $\left( \frac{2+x}{1+x},\frac{{{x}^{2}}}{1+x} \right)$ and the odd columns are the columns of the matrix $\left( \frac{x}{1+x},\frac{{{x}^{2}}}{1+x} \right)$:
$$A=\text{ }\left( \begin{matrix}
   2 & 0 & 0 & 0 & 0 & 0 & 0 & \cdots   \\
   -1 & 1 & 0 & 0 & 0 & 0 & 0 & \cdots   \\
   1 & -1 & 2 & 0 & 0 & 0 & 0 & \cdots   \\
   -1 & 1 & -3 & 1 & 0 & 0 & 0 & \cdots   \\
   1 & -1 & 4 & -2 & 2 & 0 & 0 & \cdots   \\
   -1 & -1 & -5 & 3 & -5 & 1 & 0 & \cdots   \\
   1 & 1 & 6 & -4 & 9 & -3 & 2 & \cdots   \\
   \vdots  & \vdots  & \vdots  & \vdots  & \vdots  & \vdots  & \vdots  & \ddots   \\
\end{matrix} \right).$$

Let $a\left( x \right)$  is a formal power series , $c\left( x \right)$ is a polynomial . Denote
$$\left( a\left( x \right)|c\left( x \right) \right)=\sum\limits_{n=0}^{\infty }{{{a}_{n}}{{c}_{n}}}.$$
{\bfseries Theorem 2.}
$$\left( A{{x}^{n}}|B{{x}^{m}} \right)=2{{\delta }_{n,m}}.$$
{\bfseries Proof.} For any $a\left( x \right)$, $c\left( x \right)$ we have
$$\left( a\left( x \right)|c\left( x \right) \right)=\left( Pa\left( x \right)|{{\left( {{P}^{-1}} \right)}^{T}}c\left( x \right) \right).$$
If $Pa\left( x \right)=a\left( -x \right)$, ${{\left( {{P}^{-1}} \right)}^{T}}c\left( x \right)=-c\left( -x \right)$ or $Pa\left( x \right)=-a\left( -x \right)$, ${{\left( {{P}^{-1}} \right)}^{T}}c\left( x \right)=c\left( -x \right)$, then $\left( a\left( x \right)|c\left( x \right) \right)=0$. Thus,
$$\left( A{{x}^{2n}}|B{{x}^{2m+1}} \right)=\left( A{{x}^{2n+1}}|B{{x}^{2m}} \right)=0.$$
Since
$$A{{x}^{2n}}=2{{\left( 1+x{{\left( \log a\left( x \right) \right)}^{\prime }},xa\left( x \right) \right)}^{-1}}{{x}^{2n}},$$
$$B{{x}^{2m}}=\left[ 2m,\to  \right]\left( 1+x{{\left( \log a\left( x \right) \right)}^{\prime }},xa\left( x \right) \right),$$
$$A{{x}^{2n+1}}={{\left( a\left( x \right),xa\left( x \right) \right)}^{-1}}{{x}^{2n+1}},$$
$$B{{x}^{2m+1}}=\left[ 2m+1,\to  \right]2\left( a\left( x \right),xa\left( x \right) \right),$$
then
$$\left( A{{x}^{2n}}|B{{x}^{2m}} \right)=\left( A{{x}^{2n+1}}|B{{x}^{2m+1}} \right)=2{{\delta }_{n,m}}.$$

Note that
$$A\frac{1}{1-x}=\frac{2\left( 1+x \right)}{1+x-{{x}^{2}}}\to 2\left( 1,\text{ }0,\text{ }1,-1,\text{ }2,-3,\text{ }5,\text{ }... \right),$$
$$B\frac{1}{1-x}=\frac{2\left( 1+x \right)}{1-x-{{x}^{2}}}\to 2\left( 1,\text{ }2,\text{ }3,\text{ }5,\text{ }... \right),$$
$$A\frac{1}{1+x}=\frac{2}{1+x-{{x}^{2}}}\to 2\left( 1,-1,\text{ }2,-3,\text{ }5,\text{ }... \right),$$  
$$B\frac{-1}{1+x}=\frac{2x}{1-x-{{x}^{2}}}\to 2\left( 0,\text{ }1,\text{ }1,\text{ }2,\text{ }3,\text{ }... \right).$$
In this connection, bases $A$, $B$ we will call Fibonacci bases. Denote $\alpha ={\left( 1+\sqrt{5} \right)}/{2}\;$. Note that
$$A\frac{1+\sqrt{5}x}{1-{{x}^{2}}}=2{{\left( 1-\frac{1}{\alpha }x \right)}^{-1}}\to 2\left( 1,\frac{1}{\alpha },\frac{1}{{{\alpha }^{2}}},\frac{1}{{{\alpha }^{3}}},... \right),$$ 
 $$B\frac{\sqrt{5}+x}{1-{{x}^{2}}}=2\alpha {{\left( 1-\alpha x \right)}^{-1}}\to 2\left( \alpha ,{{\alpha }^{2}},{{\alpha }^{3}},... \right),$$
$$A\frac{1-\sqrt{5}x}{1-{{x}^{2}}}=2{{\left( 1+\alpha x \right)}^{-1}}\to 2\left( 1,-\alpha ,{{\alpha }^{2}},-{{\alpha }^{3}},... \right),$$  
$$B\frac{-\sqrt{5}+x}{1-{{x}^{2}}}=-\frac{2}{\alpha }{{\left( 1+\frac{1}{\alpha }x \right)}^{-1}}=2\left( -\frac{1}{\alpha },\frac{1}{{{\alpha }^{2}}},-\frac{1}{{{\alpha }^{3}}},... \right).$$

Denote 
$${{D}_{n}}\left( x,-1 \right)={{L}_{n}}\left( x \right),  \qquad{{E}_{n-1}}\left( x,-1 \right)={{F}_{n}}\left( x \right), \qquad{{F}_{0}}\left( x \right)=0;$$
$${{L}_{-n}}\left( x \right)={{L}_{n}}\left( -x \right)={{\left( -1 \right)}^{n}}{{L}_{n}}\left( x \right),  \qquad{{F}_{-n}}\left( x \right)={{F}_{n}}\left( -x \right)={{\left( -1 \right)}^{n-1}}{{F}_{n}}\left( x \right).$$ 
{\bfseries Theorem 3.}
$$\left[ n,\to  \right]A={{x}^{n}}\left( {{L}_{-n}}\left( \frac{1}{x} \right)+{{F}_{-n}}\left( \frac{1}{x} \right) \right),  \qquad\left[ n,\to  \right]B={{x}^{n}}\left( {{L}_{n+1}}\left( x \right)+{{F}_{n+1}}\left( x \right) \right).$$
{\bfseries Proof.}
$$\left( \frac{2+x}{1+x},\frac{{{x}^{2}}}{1+x} \right)\frac{1}{1-{{\varphi }^{2}}x}=\frac{2+x}{1+x-{{\varphi }^{2}}{{x}^{2}}}=\sum\limits_{n=0}^{\infty }{{{D}_{n}}\left( -1,-{{\varphi }^{2}} \right){{x}^{n}}}=$$
$$=\sum\limits_{n=0}^{\infty }{{{\varphi }^{n}}}{{\left( -1 \right)}^{n}}{{D}_{n}}\left( \frac{1}{\varphi },-1 \right){{x}^{n}},$$
$$\left( \frac{x}{1+x},\frac{{{x}^{2}}}{1+x} \right)\frac{\varphi }{1-{{\varphi }^{2}}x}=\frac{\varphi x}{1+x-{{\varphi }^{2}}{{x}^{2}}}=\sum\limits_{n=1}^{\infty }{\varphi }{{E}_{n-1}}\left( -1,-{{\varphi }^{2}} \right){{x}^{n}}=$$
$$=\sum\limits_{n=1}^{\infty }{{{\varphi }^{n}}{{\left( -1 \right)}^{n-1}}}{{E}_{n-1}}\left( \frac{1}{\varphi },-1 \right){{x}^{n}};$$
$$\left( 1,x\left( 1+x \right) \right)\frac{1}{1-{{\varphi }^{2}}x}=\frac{1}{1-{{\varphi }^{2}}x-{{\varphi }^{2}}{{x}^{2}}}=\sum\limits_{n=0}^{\infty }{{{E}_{n}}}\left( {{\varphi }^{2}},-{{\varphi }^{2}} \right){{x}^{n}}=$$
$$=\sum\limits_{n=0}^{\infty }{{{\varphi }^{n}}{{E}_{n}}\left( \varphi ,-1 \right){{x}^{n}}}.$$
$$\left( 1+2x,x\left( 1+x \right) \right)\frac{\varphi }{1-{{\varphi }^{2}}x}=\frac{\varphi \left( 1+2x \right)}{1-{{\varphi }^{2}}x-{{\varphi }^{2}}{{x}^{2}}}=\frac{1}{\varphi x}\left( \frac{2-{{\varphi }^{2}}x}{1-{{\varphi }^{2}}x-{{\varphi }^{2}}{{x}^{2}}}-2 \right)=$$
$$=\sum\limits_{n=0}^{\infty }{{{\varphi }^{-1}}}{{D}_{n+1}}\left( {{\varphi }^{2}},-{{\varphi }^{2}} \right){{x}^{n}}=\sum\limits_{n=0}^{\infty }{{{\varphi }^{n}}}{{D}_{n+1}}\left( \varphi ,-1 \right){{x}^{n}}.$$

From the identity $2{{A}^{-1}}={{B}^{T}}$ can be seen that the space formed by the $2n$  first columns of the matrix $A$ contains the space of polynomials of degree $<n$. Basis $B$ can be considered as a basis of the space of formal power series, but with one nuance. Denote
$$a\left( x \right)=\sum\limits_{n=0}^{\infty }{{{a}_{2n}}{{x}^{2n}}}+\sum\limits_{n=0}^{\infty }{{{a}_{2n+1}}{{x}^{2n+1}}}={{a}_{1}}\left( {{x}^{2}} \right)+x{{a}_{2}}\left( {{x}^{2}} \right).$$
Then
$$Ba\left( x \right)={{a}_{1}}\left( x+{{x}^{2}} \right)+\left( 1+2x \right){{a}_{2}}\left( x+{{x}^{2}} \right).$$
If
$$a\left( x \right)=c\left( {{x}^{2}} \right)\left( \sqrt{1+4{{x}^{2}}}-x \right),$$ 
where $c\left( x \right)$ is an arbitrary series, then $Ba\left( x \right)=0$. Thus, transformation $B$ annuls the space of  series of the form $c\left( {{x}^{2}} \right)\left( \sqrt{1+4{{x}^{2}}}-x \right)$. Since
$${{\left( 1,x\left( 1+x \right) \right)}^{-1}}=\left( 1,\frac{\sqrt{1+4x}-1}{2} \right),$$ 
$${{\left( 1+2x,x\left( 1+x \right) \right)}^{-1}}=\left( \frac{1}{\sqrt{1+4x}},\frac{\sqrt{1+4x}-1}{2} \right),$$
two Riordan matrices are right inverse matrices of the matrix $B$:
$$B_{1}^{-1}=\left( 1,\frac{\sqrt{1+4{{x}^{2}}}-1}{2} \right)=\left( \begin{matrix}
   1 & 0 & 0 & 0 & 0 & \cdots   \\
   0 & 0 & 0 & 0 & 0 & \cdots   \\
   0 & 1 & 0 & 0 & 0 & \cdots   \\
   0 & 0 & 0 & 0 & 0 & \cdots   \\
   0 & -1 & 1 & 0 & 0 & \cdots   \\
   0 & 0 & 0 & 0 & 0 & \cdots   \\
   0 & 2 & -2 & 1 & 0 & \cdots   \\
   0 & 0 & 0 & 0 & 0 & \cdots   \\
   0 & -5 & 5 & -3 & 1 & \cdots   \\
   \vdots  & \vdots  & \vdots  & \vdots  & \vdots  & \ddots   \\
\end{matrix} \right),$$
$$B_{2}^{-1}=\left( \frac{x}{\sqrt{1+4{{x}^{2}}}},\frac{\sqrt{1+4{{x}^{2}}}-1}{2} \right)=\left( \begin{matrix}
   0 & 0 & 0 & 0 & \cdots   \\
   1 & 0 & 0 & 0 & \cdots   \\
   0 & 0 & 0 & 0 & \cdots   \\
   -2 & 1 & 0 & 0 & \cdots   \\
   0 & 0 & 0 & 0 & \cdots   \\
   6 & -3 & 1 & 0 & \cdots   \\
   0 & 0 & 0 & 0 & \cdots   \\
   -20 & 10 & -4 & 1 & \cdots   \\
   \vdots  & \vdots  & \vdots  & \vdots  & \ddots   \\
\end{matrix} \right).$$
These matrices can be used to find the coordinates of an arbitrary series $a\left( x \right)$ in the basis $B$. Obviously, if the coefficients of the series $b\left( x \right)$ are the coordinates of the series $a\left( x \right)$, then the coefficients of the series $b\left( x \right)+c\left( x \right)$, $Bc\left( x \right)=0$, also are the coordinates of the series $a\left( x \right).$\\
{\bfseries Example 2.}
$${{\left( 1+x \right)}^{n}}=B{{\left( \frac{1+\sqrt{1+4{{x}^{2}}}}{2} \right)}^{n}}.$$
On the other hand, since $PA=MAM$, ${{A}^{T}}=2{{B}^{-1}}$, then
$${{\left( 1+x \right)}^{n}}=B\frac{{{x}^{n}}{{L}_{n}}\left( {1}/{x}\; \right)+{{x}^{n}}{{F}_{n}}\left( {1}/{x}\; \right)}{2}.$$
Really,
$${{\left( \frac{1+\sqrt{1+4{{x}^{2}}}}{2} \right)}^{n}}=\frac{{{x}^{n}}{{L}_{n}}\left( {1}/{x}\; \right)+{{x}^{n}}{{F}_{n}}\left( {1}/{x}\; \right)}{2}+\frac{{{x}^{n-1}}{{F}_{n}}\left( {1}/{x}\; \right)\left( \sqrt{1+4{{x}^{2}}}-x \right)}{2}.$$
\section{Generalized Fibonacci bases}
Let, ${{D}_{1}}$, ${{D}_{2}}$ are the diagonal matrices such that
$${{D}_{1}}{{x}^{2n}}=\frac{1}{2}{{x}^{2n}}, \quad{{D}_{1}}{{x}^{2n+1}}={{x}^{2n+1}},  \quad{{D}_{2}}{{x}^{2n}}={{x}^{2n}}, \quad{{D}_{2}}{{x}^{2n+1}}=\frac{1}{2}{{x}^{2n+1}}.$$
Denote $A{{D}_{1}}=\tilde{A}$, $B{{D}_{2}}=\tilde{B}$. Then
$$\tilde{A}\frac{1}{1-x}=\frac{1}{2}\left( \frac{2+3x}{1+x-{{x}^{2}}} \right)\to \frac{1}{2}\left( 2,1,1,0,1,-1,2,-3,... \right),$$
$$\tilde{B}\frac{1}{1-x}=\frac{1}{2}\left( \frac{3+2x}{1-x-{{x}^{2}}} \right)\to \frac{1}{2}\left( 3,5,8,13,21,... \right),$$
$$\tilde{A}\frac{1}{1+x}=\frac{1}{2}\left( \frac{2-x}{1+x-{{x}^{2}}} \right)\to \frac{1}{2}\left( 2,-3,5,-8,13,... \right),$$  
$$\tilde{B}\frac{-1}{1+x}=\frac{1}{2}\left( \frac{-1+2x}{1-x-{{x}^{2}}} \right)\to \frac{1}{2}\left( -1,1,0,1,1,2,3,5,... \right).$$
Bases $\tilde{A}$, $\tilde{B}$ we will call reduced Fibonacci bases. Generalization of these bases are the bases ${{A}_{\left( \varphi ,\beta  \right)}}$, ${{B}_{\left( \varphi ,\beta  \right)}}$ which are constructed as follows. Let
$$b\left( xa\left( x \right) \right)=a\left( x \right),    \qquad a\left( x{{b}^{-1}}\left( x \right) \right)=b\left( x \right),$$
where
$$a\left( x \right)=\frac{\left( {\varphi }/{2}\; \right)x+\sqrt{1+\left( {{\left( {\varphi }/{2}\; \right)}^{2}}-\beta  \right){{x}^{2}}}}{1-\beta {{x}^{2}}},$$
$$b\left( x \right)=\sqrt{1+\varphi x+\beta {{x}^{2}}},  \qquad1-x{{\left( \log b\left( x \right) \right)}^{\prime }}=\frac{1+\left( {\varphi }/{2}\; \right)x}{1+\varphi x+\beta {{x}^{2}}}.$$
Then
$${{B}_{\left( \varphi ,\beta  \right)}}{{x}^{2n}}=\left[ 2n,\to  \right]\left( 1+x{{\left( \log a\left( x \right) \right)}^{\prime }},xa\left( x \right) \right)={{\left( \beta +\varphi x+{{x}^{2}} \right)}^{n}},$$
$${{B}_{\left( \varphi ,\beta  \right)}}{{x}^{2n+1}}=\left[ 2n+1,\to  \right]\left( a\left( x \right),xa\left( x \right) \right)=\left( \frac{\varphi }{2}+x \right){{\left( \beta +\varphi x+{{x}^{2}} \right)}^{n}},$$
$${{A}_{\left( \varphi ,\beta  \right)}}{{x}^{2n}}=\left( \frac{1+\left( {\varphi }/{2}\; \right)x}{1+\varphi x+\beta {{x}^{2}}},\frac{x}{\sqrt{1+\varphi x+\beta {{x}^{2}}}} \right){{x}^{2n}}=\frac{\left( 1+\left( {\varphi }/{2}\; \right)x \right){{x}^{2n}}}{{{\left( 1+\varphi x+\beta {{x}^{2}} \right)}^{n+1}}},$$
$${{A}_{\left( \varphi ,\beta  \right)}}{{x}^{2n+1}}=\left( \frac{1}{\sqrt{1+\varphi x+\beta {{x}^{2}}}},\frac{x}{\sqrt{1+\varphi x+\beta {{x}^{2}}}} \right){{x}^{2n+1}}=\frac{{{x}^{2n+1}}}{{{\left( 1+\varphi x+\beta {{x}^{2}} \right)}^{n+1}}},$$
$$\left( {{A}_{\left( \varphi ,\beta  \right)}}{{x}^{n}}|{{B}_{\left( \varphi ,\beta  \right)}}{{x}^{m}} \right)={{\delta }_{n,m}},$$
$${{P}^{\varphi }}{{A}_{\left( \varphi ,\beta  \right)}}=M{{A}_{\left( \varphi ,\beta  \right)}}M, \qquad{{\left( {{P}^{-\varphi }} \right)}^{T}}{{B}_{\left( \varphi ,\beta  \right)}}=M{{B}_{\left( \varphi ,\beta  \right)}}M.$$

Thus, $\tilde{A}={{A}_{\left( 1,0 \right)}}$, $\tilde{B}={{B}_{\left( 1,0 \right)}}$. Note that the set of matrices ${{A}_{\left( \varphi ,\beta  \right)}}$ (respectively the set of matrices ${{B}_{\left( \varphi ,\beta  \right)}}$) contains two matrix groups. Firstly, these are the powers of Pascal matrix: ${{P}^{\varphi }}={{A}_{\left( -2\varphi ,{{\varphi }^{2}} \right)}}$, ${{\left( {{P}^{\varphi }} \right)}^{T}}={{B}_{\left( 2\varphi ,{{\varphi }^{2}} \right)}}$. Secondly, ${{A}_{\left( 0,\beta  \right)}}$, ${{B}_{\left( 0,\beta  \right)}}$:
$${{A}_{\left( 0,\beta  \right)}}{{x}^{2n}}=\frac{{{x}^{2n}}}{{{\left( 1+\beta {{x}^{2}} \right)}^{n+1}}},  \qquad{{A}_{\left( 0,\beta  \right)}}{{x}^{2n+1}}=\frac{{{x}^{2n+1}}}{{{\left( 1+\beta {{x}^{2}} \right)}^{n+1}}},$$
$${{B}_{\left( 0,\beta  \right)}}{{x}^{2n}}={{\left( \beta +{{x}^{2}} \right)}^{n}}, \qquad{{B}_{\left( 0,\beta  \right)}}{{x}^{2n+1}}=x{{\left( \beta +{{x}^{2}} \right)}^{n}},$$
$${{\left( {{A}_{\left( 0,\beta  \right)}} \right)}^{T}}={{B}_{\left( 0,-\beta  \right)}}.$$
{\bfseries Theorem 4.}
$${{A}_{\left( \varphi ,{{\beta }_{1}} \right)}}{{A}_{\left( 0,{{\beta }_{2}} \right)}}={{A}_{\left( \varphi ,{{\beta }_{1}}+{{\beta }_{2}} \right)}},  \qquad{{B}_{\left( \varphi ,{{\beta }_{1}} \right)}}{{B}_{\left( 0,{{\beta }_{2}} \right)}}={{B}_{\left( \varphi ,{{\beta }_{1}}+{{\beta }_{2}} \right)}}.$$
{\bfseries Proof.}
$$\left( \frac{1+\left( {\varphi }/{2}\; \right)x}{1+\varphi x+{{\beta }_{1}}{{x}^{2}}},\frac{{{x}^{2}}}{1+\varphi x+{{\beta }_{1}}{{x}^{2}}} \right)\left( \frac{1}{1+{{\beta }_{2}}x},\frac{x}{1+{{\beta }_{2}}x} \right)=$$
$$=\left( \frac{1+\left( {\varphi }/{2}\; \right)x}{1+\varphi x+\left( {{\beta }_{1}}+{{\beta }_{2}} \right){{x}^{2}}},\frac{{{x}^{2}}}{1+\varphi x+\left( {{\beta }_{1}}+{{\beta }_{2}} \right){{x}^{2}}} \right),$$
$$\left( \frac{x}{1+\varphi x+{{\beta }_{1}}{{x}^{2}}},\frac{{{x}^{2}}}{1+\varphi x+{{\beta }_{1}}{{x}^{2}}} \right)\left( \frac{1}{1+{{\beta }_{2}}x},\frac{x}{1+{{\beta }_{2}}x} \right)=$$
$$=\left( \frac{x}{1+\varphi x+\left( {{\beta }_{1}}+{{\beta }_{2}} \right){{x}^{2}}},\frac{{{x}^{2}}}{1+\varphi x+\left( {{\beta }_{1}}+{{\beta }_{2}} \right){{x}^{2}}} \right);$$
$$\left( 1,{{\beta }_{1}}+\varphi x+{{x}^{2}} \right)\left( 1,{{\beta }_{2}}+x \right)=\left( 1,\left( {{\beta }_{1}}+{{\beta }_{2}} \right)+\varphi x+{{x}^{2}} \right),$$
$$\left( \frac{\varphi }{2}+x,{{\beta }_{1}}+\varphi x+{{x}^{2}} \right)\left( 1,{{\beta }_{2}}+x \right)=\left( \frac{\varphi }{2}+x,\left( {{\beta }_{1}}+{{\beta }_{2}} \right)+\varphi x+{{x}^{2}} \right),$$
where we used the “generalized” Riordan matrices $\left( f\left( x \right),g\left( x \right) \right)$ for which the condition ${{g}_{0}}\ne 0$ is valid [17], [18].

Thus,  ${{A}_{\left( \varphi ,\beta  \right)}}={{A}_{\left( \varphi ,0 \right)}}{{A}_{\left( 0,\beta  \right)}}$,  ${{B}_{\left( \varphi ,\beta  \right)}}={{B}_{\left( \varphi ,0 \right)}}{{B}_{\left( 0,\beta  \right)}}$. Since
$$\left[ n,\to  \right]{{A}_{\left( \varphi ,0 \right)}}={{x}^{n}}\left( \frac{1}{2}{{L}_{-n}}\left( \frac{\varphi }{x} \right)+{{F}_{-n}}\left( \frac{\varphi }{x} \right) \right),$$  
$$\left[ n,\to  \right]{{B}_{\left( \varphi ,0 \right)}}={{x}^{n}}\left( \frac{1}{2}{{L}_{n+1}}\left( \varphi x \right)+{{F}_{n+1}}\left( \varphi x \right) \right),$$
then
$$\left[ n,\to  \right]{{A}_{\left( \varphi ,\beta  \right)}}=\frac{1}{2}\left( 1,\sqrt{{{x}^{2}}-\beta } \right){{x}^{n}}{{L}_{-n}}\left( \frac{\varphi }{x} \right)+\left( \frac{x}{\sqrt{{{x}^{2}}-\beta }},\sqrt{{{x}^{2}}-\beta } \right){{x}^{n}}{{F}_{-n}}\left( \frac{\varphi }{x} \right)=$$
$$=\frac{1}{2}{{\left( \sqrt{{{x}^{2}}-\beta } \right)}^{n}}{{L}_{-n}}\left( \frac{\varphi }{\sqrt{{{x}^{2}}-\beta }} \right)+x{{\left( \sqrt{{{x}^{2}}-\beta } \right)}^{n-1}}{{F}_{-n}}\left( \frac{\varphi }{\sqrt{{{x}^{2}}-\beta }} \right);$$ 
$$\left[ n,\to  \right]{{B}_{\left( \varphi ,\beta  \right)}}=$$
$$=\frac{1}{2}\left( \frac{1}{\sqrt{1-\beta {{x}^{2}}}},\frac{x}{\sqrt{1-\beta {{x}^{2}}}} \right){{x}^{n}}{{L}_{n+1}}\left( \varphi x \right)+\left( \frac{1}{1-\beta {{x}^{2}}},\frac{x}{\sqrt{1-\beta {{x}^{2}}}} \right){{x}^{n}}{{F}_{n+1}}\left( \varphi x \right)=$$
$$=\frac{{{x}^{n}}}{2{{\left( \sqrt{1-\beta {{x}^{2}}} \right)}^{n+1}}}{{L}_{n+1}}\left( \frac{\varphi x}{\sqrt{1-\beta {{x}^{2}}}} \right)+\frac{{{x}^{n}}}{{{\left( \sqrt{1-\beta {{x}^{2}}} \right)}^{n+2}}}{{F}_{n+1}}\left( \frac{\varphi x}{\sqrt{1-\beta {{x}^{2}}}} \right).$$
{\bfseries Example 3.} Since ${{P}^{\varphi }}={{A}_{\left( -2\varphi ,{{\varphi }^{2}} \right)}}$, ${{\left( {{P}^{\varphi }} \right)}^{T}}={{B}_{\left( 2\varphi ,{{\varphi }^{2}} \right)}}$, then
$${{\left( \varphi +x \right)}^{n}}=$$
$$=\frac{1}{2}{{\left( \sqrt{{{x}^{2}}-{{\varphi }^{2}}} \right)}^{n}}{{L}_{n}}\left( \frac{2\varphi }{\sqrt{{{x}^{2}}-{{\varphi }^{2}}}} \right)+x{{\left( \sqrt{{{x}^{2}}-{{\varphi }^{2}}} \right)}^{n-1}}{{F}_{n}}\left( \frac{2\varphi }{\sqrt{{{x}^{2}}-{{\varphi }^{2}}}} \right),$$
$$\frac{1}{{{\left( 1-\varphi x \right)}^{n+1}}}=$$
$$=\frac{1}{2{{\left( \sqrt{1-{{\varphi }^{2}}{{x}^{2}}} \right)}^{n+1}}}{{L}_{n+1}}\left( \frac{2\varphi x}{\sqrt{1-{{\varphi }^{2}}{{x}^{2}}}} \right)+\frac{1}{{{\left( \sqrt{1-{{\varphi }^{2}}{{x}^{2}}} \right)}^{n+2}}}{{F}_{n+1}}\left( \frac{2\varphi x}{\sqrt{1-{{\varphi }^{2}}{{x}^{2}}}} \right).$$

E-mail: {evgeniy\symbol{"5F}burlachenko@list.ru}
\end{document}